\documentclass[letterpaper,12pt]{article}
\usepackage{amssymb}
\usepackage{amsmath}
\usepackage{amscd}
\usepackage[all]{xy}

\newtheorem{theo}{Theorem}[section]
\newtheorem{prop}[theo]{Proposition}
\newtheorem{lem}[theo]{Lemma}
\newtheorem{cor}[theo]{Corollary}
\newtheorem{defi}[theo]{Definition}

\def \Br {{\rm{Br}}}

\def \Ga {{\Gamma}}

\def \Pic {{\rm {Pic}}}

\def \Gal {{\rm{Gal}}}
\def \Ker {{\rm{Ker\,}}}
\def \Im {{\rm {Im\,}}}
\def \Fr{{\rm{Fr}}}

\def \A{{\mathbb A}}
\def \P{{\mathbb P}}
\def \Spec {{\rm{Spec\,}}}
\def \dim {{\rm{dim\,}}}

\def \End {{\rm {End}}}
\def \Pic {{\rm {Pic}}}
\def \GL {{\rm {GL}}}

\def \Aut{{\rm Aut}}

\def\ov{\overline}

\def \Z {{\mathbb Z}}
\def \Q {{\mathbb Q}}
\def \F {{\mathbb F}}

\def \rk {{\rm{rk}}}

\def\G{{\mathbb G}}

\def\C{{\mathbb C}}

\def\lra{\longrightarrow}

\def\H{{\rm H}}

\def\O{{\cal O}}

\def\NS{{\rm NS}}
\def\Kum{{\rm Kum}}

\def\O{{\cal O}}

\def\res{{\rm res}}
\def\cores{{\rm cores}}

\def\ga{\gamma}

\def\et{\rm{\acute et}}

\newcommand{\bthe}{\begin{theo}}
\newcommand{\ble}{\begin{lem}}
\newcommand{\bpr}{\begin{prop}}
\newcommand{\bco}{\begin{cor}}
\newcommand{\bde}{\begin{defi}}
\newcommand{\ethe}{\end{theo}}
\newcommand{\ele}{\end{lem}}
\newcommand{\epr}{\end{prop}}
\newcommand{\eco}{\end{cor}}
\newcommand{\ede}{\end{defi}}

\def\ga{{\alpha}}
\def\gb{{\beta}}
\def\gg{{\gamma}}
\def\gd{{\delta}}
\def\ge{{\epsilon}}
\def\gL{{\Lambda}}
\def\gl{{\lambda}}

\def\sP{{\mathcal P}}

\def\beq{\begin{equation} \label}

\def\half{{\textstyle{\frac{1}{2}}}}

\title{On the Brauer group of diagonal quartic surfaces}
\author{Evis Ieronymou, Alexei N. Skorobogatov and Yuri G. Zarhin}
\date{(with an appendix by Sir Peter Swinnerton-Dyer)}
\begin{document}
\baselineskip=15pt
\maketitle

\begin{abstract}
We obtain an easy sufficient condition for the Brauer group of a
diagonal quartic surface $D$ over $\Q$ to be algebraic. We also give an
upper bound for the order of the quotient of the Brauer group
of $D$ by the image of the Brauer group of $\Q$. The proof is based
on the isomorphism of the Fermat quartic surface with a Kummer
surface due to Masumi Mizukami.
\end{abstract}

\noindent{Mathematics Subject Classification: 11D25, 14F22, 14J28}

\section*{Introduction}

Let $D\subset\P^3_\Q$ be the quartic surface defined by the equation
\begin{equation}
x_0^4+a_1x_1^4+a_2x_2^4+a_3x_3^4=0, \label{D}
\end{equation}
where $a_1,\,a_2,\,a_3\in\Q^*$. Let $H_D\subset\Q^*$ be the subgroup
generated by $-1$, $4$, $a_1$, $a_2$, $a_3$ and the 4-th powers $\Q^{*4}$.
Write $\ov D$ for the surface over an algebraic closure $\ov \Q$
obtained from $D$ by extending the ground field to $\ov \Q$,
and let $\Br_1(D)=\Ker[\Br(D)\to \Br(\ov D)]$.
Our first main result (Corollary \ref{main}) 
states that if $\{2,3,5\}\cap H_D=\emptyset$, then
$\Br(D)=\Br_1(D)$, that is,
the Brauer group of $D$ has no transcendental elements. Note that $\Br(D)$
is known to be finite modulo $\Br_0(D)=\Im[\Br(\Q)\to\Br(D)]$ by a
general theorem proved in \cite{SZ1}.
The complete list of possible values of the finite abelian group
$\Br_1(D)/\Br_0(D)$ can be found in the thesis of Martin Bright \cite{bright}.

Our proof is based on the crucial observation that the Fermat quartic surface
$X\subset \P^3_\Q$ given by
\begin{equation}
x_0^4+x_1^4+x_2^4+x_3^4=0 \label{X}
\end{equation}
is a Kummer surface, at least after an appropriate extension of the ground
field. Over $\C$ this was first observed with some surprise in 1971 by
I.R. Shafarevich and I.I. Piatetskii-Shapiro as an application of their
global Torelli theorem for complex K3 surfaces \cite{SPS}.
In his thesis \cite{Miz1} (see also \cite{Miz2}) Masumi Mizukami constructed
an explicit isomorphism between $X$ and the Kummer surface $\Kum(A)$
associated with a certain abelian surface $A$ over $\Q$.
The details of Mizukami's construction
can be found in the appendix to this paper written by
Peter Swinnerton-Dyer.
There is a rational isogeny $A\to E\times E$ of degree 2,
where $E$ is the elliptic curve $y^2=x^3-4x$.
The Kummer surface $\Kum(A)$ can be given by equations
(\ref{Kum(A)}) of the Appendix. Note that
Mizukami's isomorphism $X\tilde\lra\Kum(A)$ is only defined 
over $\Q(\sqrt{-1},\sqrt{2})=\Q(\mu_8)$.  
Using \cite[Prop. 1.4]{SZ2} we conclude that
the Brauer groups $\Br(\ov A)$ and $\Br(\ov X)$ are isomorphic as
modules under the Galois group $\Gal(\ov \Q/\Q(\mu_8))$.
This allows us to control torsion of odd order in $\Br(D)$, see
Theorem \ref{diag} below.
The 2-primary torsion subgroup of $\Br(D)$ was studied in
the thesis of the first named author. The result that concerns us here
is \cite[Thm. 5.2]{Evis} which states that if $2\notin H_D$,
then the $2$-primary subgroup of $\Br(D)/\Br_1(D)$ is zero.
This gives Corollary \ref{main}.

Let us note in this connection that Martin Bright listed
many diagonal quartics $D$ over $\Q$ that are everywhere locally
soluble but have no rational point of height less than $10^4$, while
$\Br_1(D)=\Br(\Q)$, see
\cite{bright}, Appendices B and C. An inspection of his tables
reveals that in all cases we have $2\in H_D$. This
hints at the possibility that a potential failure of the Hasse
principle may be explained by the Brauer--Manin obstruction
attached to a transcendental element of $\Br(D)$.
See \cite[Example 3.3]{bright2} for another example of an
everywhere locally soluble diagonal quartic with no 
algebraic Brauer--Manin obstruction
and no known rational points.

As an application of Corollary \ref{main} we exhibit diagonal quartic
surfaces $D$ over $\Q$ such that $\Br(D)=\Br(\Q)$.
Indeed, Bright's computations
(\cite{bright}, Appendix A, case A161 and its subcases) show that
$\Br_1(D)=\Br(\Q)$ for the following diagonal quartics $D$:
\begin{equation}
x_0^4+4x_1^4+cx_2^4-cx_3^4=0. \label{triv}
\end{equation}
By combining this with our Corollary \ref{main}
we see that $\Br(D)=\Br(\Q)$ for $c=1,\,6,\,7,\,9,\,10,\,11\,,\ldots$
The surfaces (\ref{triv}) have obvious $\Q$-points $(0:0:1:1)$, and
it is an interesting question
whether weak approximation holds for these surfaces.

An analysis of the Galois representations on points of order
$3$, $5$ and $16$ of the lemniscatic elliptic curve $E$, 
together with Mizukami's isomorphism
and \cite[Prop. 1.4]{SZ2}, allows one to obtain an upper bound on the
size of the Brauer group of $D$. The second main result of this paper,
Corollary \ref{c3}, says that $\Br(D)/\Br_1(D)\subset (\Z/n)^2$, where
$n=2^{10}\cdot 3\cdot 5$.
Combining this with Bright's computations \cite{bright} we obtain that
the order of $\Br(D)/\Br_0(D)$ divides $2^{25}\cdot 3^2\cdot 5^2$.
By a recent theorem of A. Kresch and Yu. Tschinkel \cite[Thm. 1]{KT}
this implies that the Brauer--Manin set $D(\A_\Q)^\Br$
is effectively computable, see Corollary \ref{c4}.

\section{Brauer group and finite morphisms}

Let $k$ be a field of characteristic $0$ with an algebraic closure
$\ov k$ and the absolute Galois group $\Ga_k=\Gal(\ov k/k)$.
If $A$ is an abelian group, we write $A_n$ for the kernel of 
the multiplication by $n$ map $A\to A$.

\bpr \label{p1}
Let $X$ and $Y$ be geometrically irreducible smooth varieties over
$k$, and let $f:Y\to X$ be a dominant, generically
finite morphism of degree $d$.
Then the kernel of the natural map $f^*:\Br(X)\to\Br(Y)$
is killed by $d$. In particular,
for any integer $n>1$ coprime to $d$
the map $f^*:\Br(X)_n\to\Br(Y)_n$ is injective.
\epr
{\em Proof} By a general theorem of Grothendieck
(see \cite{EC}, Example III.2.22, p. 107)
the embedding of the generic point
$\Spec(k(X))$ in $X$ induces an injective map $\Br(X)\hookrightarrow
\Br(k(X))$, and similarly for $Y$. Since the composition of
restriction and corestriction
$$\cores_{k(Y)/k(X)}\circ\res_{k(Y)/k(X)}:\Br(k(X))\to\Br(k(Y))\to\Br(k(X))$$
is the multiplication by $d$,
the kernel of the natural map $f^*:\Br(X)\to\Br(Y)$
is killed by $d$, so our statement follows. QED

\bco \label{c1}
A degree $d$ isogeny of abelian varieties $f:A_1\to A_2$ induces a
surjective map
of $\Ga_k$-modules $f^*:\Br(\ov A_2)\to\Br(\ov A_1)$
such that $d\,\Ker(f^*)=0$. In particular, this map induces
an isomorphism on the subgroups of elements of order
coprime to $d$.
\eco
{\em Proof} If $\ov A$ is an abelian variety over $\ov k$, 
then the N\'eron--Severi group
$\NS(\ov A)$ is torsion free. Let $r=\dim \ov A$, and let
$\rho=\rk\,\NS(\ov A)$.
Then we have
$\Br(\ov A)\simeq (\Q/\Z)^m$, where $m=r(2r-1)-\rho$, 
and $r(2r-1)$ is the second Betti number of $\ov A$,
see \cite[II]{Gr}, Cor. 3.4 on p. 82, and formula (8.9) in 
\cite[III]{Gr} on p. 146.
By Proposition \ref{p1} the map $f^*:\Br(\ov A_2)\to\Br(\ov A_1)$
is a homomorphism $(\Q/\Z)^m\to (\Q/\Z)^m$ whose kernel
is killed by $d$. Such a homomorphism is necessarily surjective,
as shows the following well known lemma. QED

\ble
Any homomorphism $(\Q/\Z)^m \to (\Q/\Z)^m$ with finite kernel is 
surjective.
\ele
{\em Proof} Let $j:(\Q/\Z)^m \to (\Q/\Z)^m$ be a homomorphism
such that $d\,\Ker(j)=0$ for a positive integer $d$. The group
$(\Q/\Z)^m$ is the union of finite subgroups 
$F_r=(\frac{1}{r}\Z/\Z)^m$ for all positive integers $r$.
We have $j(F_{d^m r})\subset F_{d^m r}$, moreover, the index of 
$j(F_{d^m r})$ in $F_{d^m r}$ divides $d^m$. This implies that
$j(F_{d^m r})$ contains  $d^m F_{d^{m} r}=F_r$.
Since this holds for all $r$, the map $j$ is surjective. QED

\bthe \label{p2}
Let $X$ and $Y$ be
geometrically irreducible smooth varieties over $k$.
Let $f:Y\to X$ be a finite flat morphism of degree $d$,
such that $k(Y)$ is a Galois extension of $k(X)$ with
Galois group $G$.
Then $d^2\Br(Y)^G\subset f^*\Br(X)$. In particular,
for any integer $n>1$ coprime to $d=|G|$
the natural map $f^*:\Br(X)_n\to\Br(Y)_n^G$ is an isomorphism.
\ethe
{\em Proof} Let $\O_X$ and $\O_Y$ be the structure sheaves.
See \cite{Mum}, Lecture 10, for a construction of 
a natural map of coherent sheaves $f_*\O_Y\to \O_X$
which induces the norm map on the generic fibres
$k(Y)\to k(X)$. The composition
of the canonical map $\O_X\to f_*\O_Y$ with 
$f_*\O_Y\to \O_X$ sends $u$ to $u^d$. The \'etale sheaf
$\G_{m,X}$ is defined by setting $\G_{m,X}(U)=\Ga(U,\O_U)^*$
for any \'etale morphism $U\to X$, and similarly for $\G_{m,Y}$.
We thus obtain natural morphisms of sheaves
$$ \G_{m,X}\to f_*\G_{m,Y}\to\G_{m,X}, $$
whose composition sends $u$ to $u^d$. Applying
$\H^2_{\et}(X,\cdot)$ we define the maps
$$\xymatrix{\Br(X)\ar[rr]^-{\res_{Y/X}}&&
\H^2_{\et}(X,f_*\G_{m,Y})\ar[rr]^-{\cores_{Y/X}}&&\Br(X),}$$
whose composition is the multiplication by $d$.
Note that $f^*:\Br(X)\to\Br(Y)$ is the composition
of $\res_{Y/X}$ and the canonical map
\begin{equation}
\H^2_{\et}(X,f_*\G_{m,Y}) \lra \H^2_{\et}(Y,\G_{m,Y}) \label{isom}
\end{equation}
from the Leray spectral sequence \cite[Thm. 1.18 (a)]{EC}
$$\H^p_{\et}(X,{\rm R}^qf_*\G_{m,Y})\Rightarrow \H^{p+q}_{\et}(Y,\G_{m,Y}).$$
We have ${\rm R}^if_*\G_{m,Y}=0$ for all $i>0$
because $f$ is a finite morphism \cite[Cor. II.3.6]{EC}. 
Thus the Leray spectral sequence shows 
that (\ref{isom}) is an isomorphism. Therefore, we obtain the maps
$$\xymatrix{\Br(X)\ar[rr]^-{\res_{Y/X}=f^*}&& \Br(Y)\ar[rr]^-{\cores_{Y/X}}&&\Br(X).}$$

As was mentioned above, the embedding of the generic point
into $X$ induces an injective map $\Br(X)\hookrightarrow
\Br(k(X))$, and a similar map for $Y$. By functoriality we get 
the following commutative diagram
\begin{equation}
\xymatrix{
\Br(X)  \ar[d]_{\res_{Y/X}} \ar@{^{(}->}[r]& \Br(k(X)) \ar[d]^{\res_{k(Y)/k(X)}}\\
\Br(Y)  \ar[d]_{\cores_{Y/X}} \ar@{^{(}->}[r]&
\Br(k(Y)) \ar[d]^{\cores_{k(Y)/k(X)}}\\
\Br(X)  \ar@{^{(}->}[r]& \Br(k(X))}
\label{dm} \end{equation}
Let $\Ga_{k(X)}=\Gal(\ov{k(X)}/k(X))$ and $\Ga_{k(Y)}=\Gal(\ov{k(X)}/k(Y))$,
so that $\Ga_{k(X)}/\Ga_{k(Y)}=G$,
and consider the Hochschild--Serre spectral sequence of Galois cohomology
$$\H^p(G,\H^q(\Ga_{k(Y)},\ov{k(X)}^*))\Rightarrow \H^{p+q}(\Ga_{k(X)},\ov{k(X)}^*).$$
By Hilbert's Theorem 90 we have $\H^1(\Ga_{k(Y)},\ov{k(X)}^*)=0$.
We thus obtain the following exact sequence
$$\Br(k(X))\to (\Br(k(Y)))^G\to \H^3(G,k(Y)^*).$$
The last term is an abelian group killed by $|G|=d$.
This implies that for any $\alpha\in \Br(Y)^G$ we have 
$d\alpha=\res_{k(Y)/k(X)}(\gamma)$
for some $\gamma\in \Br(k(X))$. Then we have
$$
d\gamma=\cores_{k(Y)/k(X)}\circ\res_{k(Y)/k(X)}(\gamma)=\cores_{k(Y)/k(X)}(d\alpha)=
d\,\cores_{Y/X}(\alpha) \in \Br(X),$$
where the last equality is due to commutativity of the lower square of (\ref{dm}).
From the commutativity of the upper square of (\ref{dm}) we finally obtain
$$d^2\alpha=d\, \res_{k(Y)/k(X)}(\gamma)=\res_{k(Y)/k(X)}(d\gamma)=\res_{Y/X}(d\gamma)
\in f^*\Br(X).$$

For the last statement, the surjectivity is clear
since $\Br(Y)^G_n\subset d^2\Br(Y)^G$. The injectivity
follows from Proposition \ref{p1}. QED

\section{On torsion points of the lemniscata} \label{two}

Let $E$ be the lemniscatic elliptic curve $y^2=x^3-x$
over $\Q$. It has complex multiplication by $\O=\Z[i]$, where
$i=\sqrt{-1}$ acts on $E$ by sending $(x,y)$ to $(-x,iy)$.
We denote by $[a+bi]$ the complex multiplication by $a+bi\in\Z[i]$.

Let $\ell$ be a prime number, $\O_\ell=\O\otimes_\Z\Z_\ell$ and let
$T_\ell(E)$ be the $\ell$-adic Tate module of $E$.
For a subfield $K\subset\ov \Q$ we write $\Ga_K=\Gal(\ov\Q/K)$. Let
$\rho_\ell:\Ga_\Q\to \Aut_{\Z_\ell}(T_\ell(E))$ be the $\ell$-adic Galois
representation attached to $E/\Q$.

The action of $\O$ on $\ov E=E\times_\Q\ov\Q$ endows $T_\ell(E)$ 
with the natural structure of an $\O_\ell$-module; it is known that  
this $\O_\ell$-module is free of rank 1 
(see \cite{ST}, Remark on p. 502). The action
of $\O$ on $\ov E$ is defined over $\Q(i)$, and we have
$$\rho_\ell(\Ga_{\Q(i)})\subset\O_\ell^*\subset \Aut_{\Z_\ell}(T_\ell(E))$$
(\cite{ST}, Cor. 2 on p. 502), in particular, 
$\rho_\ell(\Ga_{\Q(i)})$ is abelian. In fact, by \cite{Serre1972},
p. 302, $\rho_\ell(\Ga_{\Q(i)})$ is an open subgroup of $\O_\ell^*$. 

A prime $p$ splits in $\O$ if and only if $p\equiv 1 \bmod 4$. 
Such a prime is uniquely written as $p=(a+bi)(a-bi)$, where
$a\pm bi\equiv 1\bmod 2+2i$. The principal ideals $(a+bi)$
and $(a-bi)$ of $\O$ are complex conjugate,
with residue fields isomorphic to $\F_{p}$.

Assume that $p \ne \ell$. Since $E$ has good reduction at $p$, the 
$\ell$-adic representation $\rho_{\ell}:\Ga_\Q\to \Aut_{\Z_\ell}(T_\ell(E))$
is unramified at $p$. A Frobenius element 
$\Fr_p\in \rho_{\ell}(\Ga_\Q)$ is the image of a Frobenius automorphism
at the prime $p$, and so $\Fr_p$ is well defined up to conjugation 
in $\rho_{\ell}(\Ga_\Q)$ (see \cite{SerreELL}, Ch. 1, Sect. 1.2 
and Sect. 2 for more details). The representation
$\rho_{\ell}:\Ga_{\Q(i)}\to \Aut_{\Z_\ell}(T_\ell(E))$ 
is unramified at $(a+bi)$ and $(a-bi)$, and the corresponding Frobenii
are well defined elements of the abelian group $\rho_{\ell}(\Ga_{\Q(i)})$.
In $\rho_{\ell}(\Ga_\Q)$ these two elements are conjugate by $\rho_\ell(c)$,
where $c\in \Ga_\Q$ is the complex conjugation, so they are precisely 
the elements of the conjugacy class of $\Fr_p$ in $\rho_{\ell}(\Ga_\Q)$.

A well known fact going back to the last entry of Gauss's 
mathematical diary (via Deuring's interpretation on Hecke characters) 
is that the Frobenius element in 
$\rho_\ell(\Ga_{\Q(i)})\subset \O_\ell^*$ attached to the prime ideal  
$(a+ bi)$ equals $a+ bi$ (and similarly for $a-bi$, see \cite{IR},
Thm. 5 on p. 307, or \cite{RS}, Prop. 4.1 and its proof, and Thm. 5.6). 
In what follows $\Fr_p$ stands for either $a+bi$ or $a-bi$,
for example $\Fr_5=-1+ 2i$ and $\Fr_{17}=1+ 4i$. 

We choose a basis of the free $\Z_\ell$-module $T_\ell(E)$
of rank 2 so that
the image of $[i]\in \O$ in $\End_{\Z_\ell}(T_\ell(E))$
is represented by the matrix
\[ \left( \begin{array}{rr}
0 & -1  \\
1 & 0  \end{array} \right)\]
Then $\O_\ell\subset \End_{\Z_\ell}(T_\ell(E))$ consists of the matrices
$$
\left( \begin{array}{rr}
a & -b  \\
b & a  \end{array} \right) $$
for all $a,\,b\in \Z_\ell$.

\bpr \label{j} Let $k$ be a Galois extension of $\Q(i)$.

{\rm (a)} If the exponent of $\Gal(k/\Q(i))$ divides $4$, then
$\End_{\Ga_k}(E_{\ell})=\O/\ell$ for any prime $\ell \geq 7$.

{\rm (b)} We have
$$\End_{\Ga_k}(E_{5})=\left \{ \begin{array}{ll}
\O/5 & \textrm{if}  \ \sqrt[4]{5} \notin k,  \\
\End(E_{5}) & \textrm{otherwise;}
\end{array}
\right .
 \quad
 \End_{\Ga_k}(E_{3})=\left \{ \begin{array}{ll}
\O/3 & \textrm{if}  \ \sqrt[4]{-3} \notin k,  \\
\End(E_{3}) & \textrm{otherwise.}
\end{array}
\right .$$
\epr
{\em Proof} Let $\ov\rho_\ell:\Ga_k\to \Aut_{\F_\ell}(E_\ell)\simeq \GL(2,\F_\ell)$
be the Galois representation modulo $\ell$ attached to $E/k$.
Define $\Lambda\subset\Ga_k$ as $\ov\rho_\ell^{\,-1}(\F_\ell^*)$,
and let $M=\ov \Q^\Lambda$.
If $M\not=k$, then there exists $\gamma\in\Ga_k$ such that $\ov\rho_\ell(\gamma)$
has two distinct eigenvalues
in $\ov\F_\ell$. The centralizer of $\ov\rho_\ell(\gamma)$ in
$\End(E_{\ell})$ is an $\F_\ell$-vector space of
dimension 2 which contains $\O/\ell$ and so is equal to it.
Hence in this case $\End_{\Ga_k}(E_{\ell})=\O/\ell$.
If $M=k$, then the image of $\Ga_k$ in $\GL(2,\F_\ell)$
is the group of scalar matrices, so that $\End_{\Ga_k}(E_{\ell})=\End(E_{\ell})$.

To prove (a) we note that $5$ splits in $\Q(i)$ and hence
$\Fr_5=-1+2i\in \O_\ell^*$ belongs to $\rho_\ell(\Ga_{\Q(i)})$.
Our assumption implies that ${\rm Fr}_5^4$
belongs to $\rho_\ell(\Ga_k)$.
Since ${\rm Fr}_5^4=-7+ 24i$ is not congruent to an element of $\F_\ell$
modulo $\ell$, we see that $\ov\rho_\ell(\Ga_k)\not\subset\F_\ell^*$
so that $\End_{\Ga_k}(E_{\ell})=\O/\ell$.

To prove (b) it suffices to show that when $k=\Q(i)$, then
$M=k(\sqrt[4]{5})$ for $\ell=5$, and $M=k(\sqrt[4]{-3})$ for $\ell=3$.

The case $\ell=5$.

Since $5$ splits in $\O$, the $\Ga_k$-module $E_5$ is the direct sum of characters
$\chi_1\oplus\chi_2$ with values in $\F_5^*$.
Then $M$ is the fixed field of $\Ker(\chi_1\chi_2^{-1})$.

The duplication formula gives the $x$-coordinate of the double of
a point $(x,y)$ on $E$ as $(x^2+1)^2/4x(x^2-1)$, see \cite{Sil},
Ch. X, Section 6, pp. 309--310. Using this it is easy to see that
a point $(x_1,y_1)$ such that $x_1^2=(1+2i)^{-1}$
generates $\Ker [1-2i]$, and that a point $(x_2,y_2)$ such that $x_2^2=(1-2i)^{-1}$
generates $\Ker [1+2i]$. This implies that $y_1^4=-4(1+2i)^{-3}$,
$y_2^4=-4(1-2i)^{-3}$. Then $M_1=k(y_1)$ and $M_2=k(y_2)$ are cyclic
extensions of $k$ of degree 4 which are linearly disjoint
since $M_1$ is totally ramified at the principal prime ideal
$(1+2i)$ and unramified at $(1-2i)$, while $M_2$ is totally ramified
at $(1-2i)$ and unramified at $(1+2i)$.
We can therefore identify $\Gal(M_1M_2/k)$
with $\Gal(M_1/k)\times \Gal(M_2/k)$.
Let $g_1$ denote the generator of
$\Gal(M_1/k)\simeq \Z/4$ such that $g_1(y_1)=i y_1$.
We define $g_2$ similarly. From the above it is clear that
$M$ is the fixed subfield of $g_1g_2^{-1}$.
Note that $\frac{5}{2}y_1y_2$ is fixed by $g_1g_2^{-1}$ and
$(\frac{5}{2}y_1y_2)^4=5$. Since [$M:k]=4$
we conclude that $M=\Q(\sqrt{-1},\sqrt[4]5)$.

The case $\ell=3$.

Since $3$ is inert in $\O$, $\Ga_k$ acts on $E_3$ by a character
$\chi$ with values in $\F_{9}^*$. Recall that
$\Lambda\subset \Ga_k$ is $\chi^{-1}(\pm 1)$, and $M=\ov \Q^\Lambda$.

Applying the duplication formula we immediately see that
if $P=(x,y)$ is a point of order 3 in $E$, then $x$
is a root of the polynomial $f(t)=t^4 - 2t^2 - 1/3$.
By Eisenstein's criterion $z^4+6z^2-3$ is irreducible over $\Q$,
and even over $k$ since $3$ is an irreducible element of the
unique factorisation domain $\Z[i]$. The polynomial $z^4+6z^2-3$
completely splits in $k(\sqrt[4]{-3})$ since it has a root
$(1+i)a(a^2-i)/2$, where $a=\sqrt[4]{-3}$, and $k(\sqrt[4]{-3})$
is a Galois extension of $k$.
Hence $f(t)$ is irreducible over $k$ with
splitting field $M_1=\Q(i,\sqrt[4]{-3})$.
Let $M_2=M_1(y)=M_1(\sqrt{x^3-x})$. Since $P$ has order 3, the points
$P$ and $[i]P$ span the $\F_3$-vector space $E_3$, so that
$E_3\subset E(M_2)$. It is clear that $[M_2:k]$ is $8$ or $4$.
The prime $17$ splits in $\Q(i)=k$, and hence $1+4i\in \O_3^*$
belongs to $\rho_3(\Ga_k)$.
Since $1+ 4i$ modulo $3$ has multiplicative
order $8$, the order of the Galois group $\Gal(k(E_3)/k)$
is divisible by $8$. Therefore, $M_2=k(E_3)$ is an extension of $k$
of degree $8$, $[M_2:M_1]=2$,
and $\Gal(k(E_3)/k)=(\Z[i]/3)^*=\F_9^*$ is a cyclic group of order $8$.

The $M_1$-linear automorphism of $M_2$ which maps  $y$ to $-y$
corresponds to the multiplication by $-1$ in $E_3$ and so belongs to $\Lambda$.
Therefore $M=\ov \Q^\Lambda\subset M_1$, and in fact $M=M_1$
since $\F_3^*$ has index 4 in $\F_9^*$. Thus $M=\Q(i,\sqrt[4]{-3})$. QED

\section{A sufficient condition for the Brauer group of $D$
to be algebraic}

We need an easy lemma from Galois theory.

\ble \label{g}
Let $b_i,\,d\in \Q^*$, and let
$F=\Q(\sqrt{-1},\sqrt[4]{b_1},\cdots,\sqrt[4]{b_n})$.
Then $t^4-d$ splits in $F$ if and only
if $d$ belongs to the subgroup of $\Q^*/\Q^{*4}$
generated by the classes of $-4$ and the $b_i$, $i=1,\ldots,n$.
\ele
{\em Proof} This is \cite[Lemma 5.4]{Evis}; we reproduce
the proof for the convenience of the reader.
The field $F$ is a 4-Kummer extension of $\Q(\sqrt{-1})$,
so $d$ is a 4-th power in $F$ if and only if $d$ belongs to
the subgroup of $\Q(\sqrt{-1})^*/\Q(\sqrt{-1})^{*4}$ generated by the
$b_i$, $i=1,\ldots,n$.
Moreover, the kernel of the natural map
\[\Q^*/\Q^{*4}\rightarrow \Q(\sqrt{-1})^*/\Q(\sqrt{-1})^{*4}\]
is a subgroup of order 2 generated by the class of $-4$. QED

\medskip

{From} now on let $k=\Q(\sqrt{-1},\sqrt{2})$ and
$F=k(\sqrt[4]{a_1},\sqrt[4]{a_2},\sqrt[4]{a_3})$, understood as
normal subfields of $\ov \Q$.

\bthe \label{diag}
Let $D\subset\P^3_{\Q}$ be the diagonal quartic surface {\rm (\ref{D})}.
Then for any prime $\ell \geq 7$ we have $\Br(\ov{D})^{\Ga_\Q}_\ell=0$.
Moreover, if  $5$ (resp. $3$) does not belong to the subgroup of
$\Q^*/\Q^{*4}$ generated by the classes of $-1,\,4,\,a_1,\,a_2,\,a_3$,
then $ \Br(\ov{D})^{\Ga_\Q}_5=0$ (resp. $\Br(\ov{D})^{\Ga_\Q}_3=0$).
\ethe
{\em Proof} Let $X$ be the surface (\ref{X}), and let $A$ be the abelian
surface defined in Theorem \ref{M}. Since
$D\times_\Q F\simeq X\times_\Q F$ the $\Ga_F$-modules
$\Br(\ov D)$ and $\Br(\ov X)$ are isomorphic.
By Mizukami's isomorphism (Theorem \ref{M}) the Fermat quartic
$X$ is isomorphic to $\Kum(A)$ over $k$, so that
$\Br(\ov X)$ and $\Br(\ov A)$ are isomorphic as $\Ga_k$-modules
\cite[Prop. 1.4]{SZ2}. Since $\ell$ is odd,
Corollary \ref{c1} now implies that
$\Br(\ov D)_{\ell^\infty}$ and $\Br(\ov E\times \ov E)_{\ell^\infty}$
are isomorphic as $\Ga_F$-modules,
so it is enough to prove that $\Br(\ov E\times \ov E)_\ell^{\Ga_F}=0$.
The $\Ga_k$-module $\H^2_{\et}(\ov E\times \ov E,\mu_\ell)$
is naturally isomorphic to $\Z/\ell\oplus\Z/\ell\oplus \End(E_\ell)$,
see e.g., \cite{SZ2}, formula (17).
The Kummer exact sequence gives rise to the well known exact sequence
of $\Ga_k$-modules
$$0\to \NS(\ov E\times\ov E)/\ell\to \H^2_{\et}(\ov E\times \ov E,\mu_\ell)
\to \Br(\ov E\times\ov E)_\ell\to 0,$$
where $\NS(\ov E\times\ov E)$ is the N\'eron--Severi group, which
is isomorphic to $\Z\oplus\Z\oplus\O$ as a $\Ga_k$-module.
The action of $\Ga_k$ on this module
is trivial because the complex multiplication on $E$ is defined over $k$.
The image of $\NS(\ov E\times\ov E)/\ell$ in
$\H^2_{\et}(\ov E\times \ov E,\mu_\ell)$ is
$\Z/\ell\oplus\Z/\ell\oplus \O/\ell$.

Note that $\ell$
is unramified in $\O$, thus $\O/\ell$ is either
$\F_{\ell}\oplus\F_{\ell}$ or the
field $\F_{\ell^2}$. In either case $\ell$
does not divide $|(\O/\ell)^{*}|$. Since
the image $G_\ell$ of $\Ga_k$ in $\Aut(E_{\ell})$
belongs to $(\O/\ell)^{*}$, we see that
$|G_{\ell}|$ is not divisible by $\ell$. 
It follows from Maschke's theorem that $E_\ell$, $\End(E_\ell)$
and $\H^2_{\et}(\ov E\times \ov E,\mu_\ell)$ are semisimple
$\Ga_k$-modules.
Therefore, we have an isomorphism of $\Ga_k$-modules
$$\End(E_\ell)\cong \O/\ell\oplus\Br(\ov E\times \ov E)_\ell,$$
where $\O/\ell$ carries trivial $\Ga_k$-action. We conclude that
$\Br(\ov E\times \ov E)^{\Ga_F}_\ell$ can be identified with
$\End_{\Ga_F}(E_\ell)/(\O/\ell)$.
Now the desired statements follow from Proposition
\ref{j} by Lemma \ref{g}. QED

\bco \label{main}
Let $H_D\subset\Q^*$ be the subgroup
generated by $-1$, $4$, $a_1$, $a_2$, $a_3$ and the $4$-th powers $\Q^{*4}$.
If $\{2,3,5\}\cap H_D=\emptyset$, then $\Br(D)=\Br_1(D)$.
\eco
{\em Proof} Since $\Br(D)$ is a torsion group,
any element $\alpha\in\Br(D)$ can be written as $\alpha=\beta+\gamma$
where $2^m\beta=0$
and $n\gamma=0$ for some $m,\,n\in \Z_{\geq 0}$, $n$ odd.
Thm. 5.2 of \cite{Evis} states that if $2\notin H_D$,
then the $2$-primary subgroup of $\Br(D)/\Br_1(D)$ is zero.
Thus our condition implies that $\beta\in\Br_1(D)$.
Also, $\gamma\in\Br_1(D)$ since $\Br(\ov D)_n^{\Ga_\Q}=0$
by Theorem \ref{diag}. QED

\section{An upper bound for $|\Br(D)/\Br_1(D)|$}

We start with the analysis of torsion of odd order in $\Br(D)/\Br_1(D)$.

\bpr\label{pp}
Let $D\subset\P^3_{\Q}$ be the diagonal quartic surface {\rm (\ref{D})}.
Then for any odd prime $\ell $ we have
$\Br(\ov D)^{\Ga_\Q}_{\ell^\infty}=\Br(\ov{D})^{\Ga_\Q}_\ell$.
\epr
{\em Proof} In the beginning of the proof of Theorem \ref{diag}
we have seen that $\Br(\ov D)_{\ell^\infty}$ and
$\Br(\ov E\times \ov E)_{\ell^\infty}$
are isomorphic as $\Ga_F$-modules. Also
in the proof of Theorem \ref{diag} we showed that
$\Br(\ov E\times \ov E)^{\Ga_F}_\ell=0$ for $\ell\geq 7$. Thus
it is enough to prove that
$\Br(\ov E\times \ov E)_{\ell^2}^{\Ga_F}=0$ is killed by $\ell$,
where $\ell=3$ or $\ell=5$.

Recall that $\O_\ell^*\subset \Aut_{\Z_\ell}(T_\ell(E))$.
Consider $\End_{\Z_\ell}(T_\ell(E))\simeq{\rm Mat}_2(\Z_\ell)$ as an $\O_\ell^*$-module
under conjugation. Since $\ell$ is odd, we can decompose this module
into a direct sum of $\O_\ell^*$-submodules
$\O_\ell\oplus\ov\O_\ell$, where
$$\ov\O_\ell=\O_\ell\left( \begin{array}{cc}
0 & 1  \\
1 & 0  \end{array} \right)=\left\{\left( \begin{array}{rr}
a & b  \\
b & -a  \end{array} \right), \ \ a,\,b\in \Z_\ell \ \right\}$$
We note that $[i]$ acts on $\ov\O_\ell$ by $-1$.
Now the exact sequence of $\Ga_{\Q(i)}$-modules
$$0\to\O_\ell/\ell^2\to\End(E_{\ell^2})\to \Br(\ov E\times \ov E)_{\ell^2}
\to 0$$
implies that the $\Ga_{\Q(i)}$-module
$\Br(\ov E\times \ov E)_{\ell^2}$ is obtained from the
$\O_\ell^*$-module $\ov\O_\ell/\ell^2$ via the map
$\Ga_{\Q(i)}\to \O_\ell^*$.

Since $17$ splits in $\O$, by Gauss's result
$1+4i\in \O_\ell^*$ is contained in $\rho_\ell(\Ga_{\Q(i)})$.
The exponent of $\Gal(F/\Q(i))$ divides 4, so we see that
$(1+4i)^4=161-240i$ belongs to $\rho_\ell(\Ga_F)$.
Let $x\in \ov\O_\ell/\ell^2$ be an element invariant under the action of
$\Ga_F$. Then $x$ commutes with $[161-240i]$.
Since $240$ is divisible by $\ell$ but not by $\ell^2$
we see that $\ell x$ is invariant under the action of $[i]$,
hence $\ell x=-\ell x$. Since $\ell$ is odd we conclude
that $\ell x=0$. QED

\medskip

To estimate 2-primary torsion in $\Br(D)/\Br_1(D)$ we need some preparations.

\ble \label{grp}
Let $G$ be a group of order $|G|=2^n$, and let $M$ be a torsion abelian
$2$-primary group which is a $G$-module.
If $M^G$ is killed by $2^m$, and $M_{4}\subset M^{G}$, then
$M$ is killed by $2^{m+n}$.
\ele
{\em Proof}
The proof is by induction on $n$. For $n=1$ let $g$ be
the non-trivial element of $G$. If $M$ contains an element $x$
of exact order $2^{m+2}$, then $2^m x$ has order $4$ and so
$2^m g(x)=2^m x$. This implies that $2^m(x+g(x))=2^{m+1}x\not=0$.
However, $x+g(x)\in M^G$ and by assumption $2^m(x+g(x))=0$
which is a contradiction.

When $n>1$, the group $G$ has
a proper normal subgroup $G_1\subset G$. Applying the
induction hypothesis two times, first to $(G/G_1, M^{G_1})$,
and then to $(G_1, M)$, we prove the induction step. QED

\bpr \label{terib}
The exponent of $\Br(\ov E  \times \ov E)_{2^{\infty}}^{\Ga_k}$
divides $8$, and that of
$\Br(\ov E \times \ov E)_{2^{\infty}}^{\Ga_F}$ divides
$2^3|\Gal(F/k)|$.
\epr
{\em Proof} The prime $17$ is congruent to 1 modulo 8,
hence it splits completely in the cyclotomic field $k=\Q(\mu_8)$.
Thus $\Fr_{17}=1+4i$ is contained in $\rho_2(\Ga_k)\subset \O_2^*$.
In our basis of $T_2(E)$ the complex multiplication
$[1+4i]$ is given by the matrix
\[
s=\left( \begin{array}{rr}
1 & -4  \\
4 & 1  \end{array} \right)
\]
For the first claim it is clearly enough to prove that for any
$\alpha\in\Br(\ov E\times\ov E)_{16}^{\Ga_F}$
we have $8\alpha=0$.
Consider the exact sequence of $\Ga_k$-modules
\[
 0\to \O/16 \to \End(E_{16}) \to \Br(\ov E  \times \ov E)_{16} \to 0.
\]
We represent $\alpha$ by a matrix
\[
A=\left( \begin{array}{ccc}
a & b  \\
c & d  \end{array} \right)\ \in\ \End(E_{16})\simeq {\rm Mat}_2(\Z/16).
\]
Then $sAs^{-1}-A\in \O/16$, so that $sA-As\in \O/16$, 
which immediately implies that $8(a-d)=8(c+b)=0$.
Thus $8A=8(a+ci)\in \O/16$, so that $8\alpha=0$.

To prove the second claim we note that $E_4\subset E(k)$. Indeed, it is well
known that for an elliptic curve $y^2=(x-c_1)(x-c_2)(x-c_3)$ 
over $\Q$ the field $\Q(E_4)$ is an extension of $\Q$ obtained by 
joining the square roots of $-1$ and $c_i-c_j$ for all $i\not=j$
(see, for example, \cite{Knapp}, Thm. 4.2 on p. 85). 
In our case $\Q(E_4)=\Q(\mu_8)=k$. This
implies that $\End(E_4)$ is a trivial $\Ga_k$-module, hence $\Br(\ov E \times
\ov E)_4= \End(E_4)/(\O/4)$ is also a trivial $\Ga_k$-module. 
Thus we can apply Lemma \ref{grp} to $G=\Gal(F/k)$ and 
$M=\Br(\ov E \times \ov E)_{2^{\infty}}^{\Ga_F}$. 
This completes the proof. QED

\medskip



\bpr \label{bound}
Let $D \subset\P^3_\Q$ be the quartic surface {\rm (\ref{D})}.
Then the exponent of $\Br (\ov D)_{2^{\infty}}^{\Ga_F}$ divides $2^{4}|\Gal(F/k)|$.
\epr
{\em Proof}
Let $X$ be the Fermat quartic surface (\ref{X}), and let $A$ be the
abelian surface defined in Theorem \ref{M}.
Because of the isomorphism $ X\times_k F\,\tilde\lra\,D\times_k F$
we can replace $D$ by $X$.
By \cite[Prop. 1.4]{SZ2} the $\Ga_F$-modules $\Br(\ov X)_{2^\infty}$
and $\Br(\ov A)_{2^\infty}$ are isomorphic. There is a degree 2 isogeny
$A\to E\times E$, so by Corollary \ref{c1} we have an exact sequence
of $\Ga_F$-modules
$$0\to (\Z/2)^n \to \Br(\ov E\times \ov E)_{2^\infty}\to
\Br(\ov A)_{2^\infty}\to 0.$$
It gives rise to the exact sequence
$$\Br(\ov E\times \ov E)_{2^\infty}^{\Ga_F}\to
\Br(\ov A)_{2^\infty}^{\Ga_F}\to \H^1({\Ga_F},(\Z/2)^n).$$
The last term is killed by $2$. On the other hand,
by Proposition \ref{terib} the exponent of
$\Br(\ov E \times \ov E)_{2^{\infty}}^{\Ga_F}$ divides $2^{3}|\Gal(F/k)|$.
Hence $\Br(\ov A)_{2^{\infty}}^{\Ga_F}$ is killed by $2^4|\Gal(F/k)|$. QED

\bco \label{c2}
Let $D \subset\P^3_\Q$ be the quartic surface {\rm (\ref{D})}.
Then the exponent of $\Br(\ov D)_{2^{\infty}}^{\Ga_\Q}$ divides $2^{10}$.
\eco
{\em Proof} In the notation of Proposition \ref{bound} the Galois
group $\Gal(F/k)$ is a quotient of $(\Z/4)^3$, and so
$\Br (\ov D)_{2^{\infty}}^{\Ga_F}$ is killed by $2^{10}$. Now the statement
follows from the obvious inclusion
$\Br (\ov D)_{2^{\infty}}^{\Ga_\Q}\subset\Br (\ov D)_{2^{\infty}}^{\Ga_F}$.
QED

\bco \label{c3}
Let $D \subset\P^3_\Q$ be the quartic surface {\rm (\ref{D})}.
Then

{\rm (i)} the exponent of the group $\Br(D)/\Br_1(D)$ 
divides $2^{10}\cdot 3\cdot 5$;

{\rm (ii)} the order of $\Br(D)/\Br_1(D)$ divides $2^{20}\cdot 3^2\cdot 5^2$;

{\rm (iii)} 
the order of $\Br(D)/\Br_0(D)$ divides $2^{25}\cdot 3^2\cdot 5^2$.
\eco
{\em Proof} (i) The case of $\ell$-primary torsion,
where $\ell$ is an odd prime, follows from Theorem \ref{diag}
combined with Proposition \ref{pp}. The case of 2-primary torsion is
dealt with in Corollary \ref{c2}.

(ii) It is well known that $\Pic(\ov D)=\NS(\ov D)\simeq\Z^{20}$,
see, e.g. \cite[Lemma 1]{PSD}. Since the second Betti number
of $\ov D$ is 22, we conclude that $\Br(\ov D)\simeq (\Q/\Z)^2$
(using \cite[II]{Gr}, Cor. 3.4 on p. 82, 
and formula (8.12) of \cite[III]{Gr} on p. 147).
Thus (ii) follows from (i).

(iii)
This statement follows from M. Bright's computations
that the order of $\Br_1(D)/\Br_0(D)$ divides $2^5$,
see \cite{bright}. QED

\medskip

We refer the reader to \cite{Sk}, Section 5.2,
for the definition of the Brauer--Manin set $D(\A_\Q)^\Br$.

\bco \label{c4}
Let $D \subset\P^3_\Q$ be the quartic surface {\rm (\ref{D})}.
Then the Brauer--Manin set $D(\A_\Q)^\Br$ is effectively computable.
\eco
{\em Proof} According to \cite[Thm. 1]{KT}
for a family of smooth projective surfaces $Z$ over $\Q$
defined by explicit equations, such that $\Pic(\ov Z)$
is torsion free and generated by finitely many explicitly given
divisors, the Brauer--Manin set $Z(\A_k)^\Br$ is effectively
computable whenever one has
a uniform bound on the order of $\Br(Z)/\Br_0(Z)$. The geometric
Picard group $\Pic(\ov D)\simeq\Z^{20}$ of a diagonal quartic surface
is generated by the obvious 48 lines on it \cite[Lemma 1]{PSD}.
Thus the statement follows from Corollary \ref{c3}. QED

\bigskip

\noindent{\it Acknowledgements.} We are very grateful to Tetsuji Shioda for
sending us a copy of Mizukami's thesis \cite{Miz1}, and to
Peter Swinnerton-Dyer for the Appendix and for pointing out
equation (\ref{triv}) to us. We thank Victor Flynn for his help with computations
and Alice Silverberg for her help with references.

\appendix

\section{The Fermat quartic as a Kummer surface (after Mizukami), \\
by Sir Peter Swinnerton-Dyer}

Let $k$ be a field of characteristic different from $2$.
Let $C$ be the elliptic curve over $k$ 
which is a smooth projective model of
the affine curve $v^2=(u^2-1)(u^2-\frac{1}{2})$. 
The base point $O$ of $C$ is that point at infinity at which $v/u^2=1$.

\bthe[M. Mizukami, 1977] \label{M}
Let $k$ be a field of characteristic not equal to $2$
that contains the $8$-th roots of unity.
Let $\tau:C\to C$ be the fixed point free involution
changing the signs of $v$ and $u$.
Let $A$ be the abelian surface obtained as the quotient
of $C\times C$ by the simultaneous action of $\tau$ on both
factors. Then there is an isomorphism $X\,\tilde\lra\,K=\Kum(A)$,
where $X$ is the Fermat quartic surface {\rm (\ref{X})}.
\ethe

The curves $C$ and $C'=C/\tau$ considered as elliptic curves over $\Q$
have Cremona labels 64a2 and 64a1, respectively; these curves have
good reduction away from $2$.
The short Weierstrass equation of $C'$ is
$y^2=x^3-4x$, so over $\Q(\mu_8)=\Q(\sqrt{-1},\sqrt{2})$
the curve $C'$ is isomorphic
to the elliptic curve $E$ with equation $y^2=x^3-x$.
Thus there is a degree 2 isogeny $A\to E\times E$ defined over $\Q(\mu_8)$.

\medskip

\noindent{\em Proof of Theorem} \ref{M}.
Let $T$ be that point at infinity at which $v/u^2=-1$. We shall need to consider
two copies of $C$; we distinguish these and
the associated variables by the subscripts 1 and 2.
The involution on $C_1\times C_2$
which reverses the signs of all four variables $u_1,v_1,u_2,v_2$
has no fixed points;
so it is a translation by $T_1\times T_2$
(which is the image of $O_1\times O_2$)
and $T_1\times T_2$ must be a 2-division point. Now
\[ A=C_1\times C_2/\{O_1\times O_2, T_1\times T_2\} \]
is an abelian surface equipped with a map $C_1\times C_2\rightarrow A$
of degree 2. Its function field is the field of functions
even in the four variables
$u_1,u_2,v_1,v_2$ collectively.
The involution $P\mapsto-P$ reverses the signs of
$u_1$ and
$u_2$; so the function field of $K$ over a field $k$ can be written as $k(w_1,
w_2,y,z)$ where
\[ w_1=u_1^2, \quad w_2=u_2^2, \quad y=\frac{v_1(w_2-1)}{v_2(w_1-\half)},
\quad z=u_2/u_1. \]
(The reason for the unnatural-looking choice for $y$ will appear at (\ref{E3}).)
Thus up to birational transformation we can take $K$ to be given by
\begin{equation} y^2=(w_1-1)(w_2-1)/(w_1-\half)(w_2-\half),
\quad z^2=w_2/w_1. \label{Kum(A)}\end{equation}
In particular $[k(K):k(w_1,w_2)]=4$.

In all that follows
we take $\ge$ to be a fixed solution of $\ge^4=-1$.

We need a notation for the lines and some of the conics on
the Fermat quartic surface $X$. (There are
actually two types of conic on $X$, but only one of them concerns us.) Let
$\mu,\nu$ be odd residue classes mod 8; then the 48 lines on $X$ can be written
\begin{align*}
 & L_{\mu\nu}:\;x_0=\ge^\mu x_1,\;x_2=\ge^\nu x_3; \\
 & M_{\mu\nu}:\;x_0=\ge^\mu x_2,\;x_1=\ge^\nu x_3; \\
 & N_{\mu\nu}:\;x_0=\ge^\mu x_3,\;x_1=\ge^\nu x_2.
\end{align*}
(The rejection of symmetry here is deliberate,
because cyclic symmetry plays no
part in what follows. Note also that the notation is different to Mizukami's,
and also to that of Segre.) We shall use $\gL$ to denote any one of the three
letters $L,M,N$. Then $\gL_{\ga\gb}$ meets $\gL_{\gg\gd}$ if and only
if $\ga=\gg$
or $\gb=\gd$ but not both. Conditions for $\gL_{\ga\gb}$ and
$\gL'_{\gg\gd}$ to
meet, where $\gL$ and $\gL'$ are different, are as follows:
\begin{align*}
 & L_{\ga\gb} \text{ meets } M_{\gg\gd} \text{ if and only if } \ga-\gb-\gg+\gd=0, \\
 & L_{\ga\gb} \text{ meets } N_{\gg\gd} \text{ if and only if } \ga+\gb-\gg+\gd=0, \\
 & M_{\ga\gb} \text{ meets } N_{\gg\gd} \text{ if and only if } \ga+\gb-\gg-\gd=0.
\end{align*}
Now suppose that $\gL_{\ga\gb}$ and $\gL'_{\gg\gd}$ are two lines which meet, where
$\gL$ and $\gL'$ are different; then the plane containing $\gL_{\ga\gb}$ and
$\gL'_{\gg\gd}$ meets $X$ residually in a conic, which we shall denote by
$[\gL_{\ga\gb}\gL'_{\gg\gd}]$. The lines $\gL_{\ga\gb}$ and $\gL'_{\gg\gd}$ meet this
conic twice. The lines which meet it once are those which meet neither $\gL_{\ga\gb}$
nor $\gL'_{\gg\gd}$. Any other conic meets the plane of $[\gL_{\ga\gb}\gL'_{\gg\gd}]$
twice, and using the previous sentences one can work out its intersection with
$[\gL_{\ga\gb}\gL'_{\gg\gd}]$.

If we write
\[ f_{\gl\mu\nu}=x_0+\ge^\gl x_1+\ge^\mu x_2+\ge^\nu x_3 \]
then the equations of $[L_{\ga\gb}M_{\gg\gd}]$ can be written
\begin{align*}
 & x_0-\ge^\ga x_1-\ge^\gg x_2+\ge^{\ga+\gd}x_3=f_{\ga+4,\gg+4,\ga+\gd}=0, \\
 & x_0^2+\ge^{\ga+\gd}x_0x_3+\ge^{2\ga+2\gd}x_3^2+\ge^{2\ga}x_1^2+\ge^{\ga+\gg}x_1x_2
+\ge^{2\gg}x_2^2=0.
\end{align*}
Thus the intersection of $X$ with $f_{\ga+4,\gg+4,\ga+\gd}=0$ is
\[ L_{\ga\gb}+M_{\gg\gd}+[L_{\ga\gb}M_{\gg\gd}] \]
where $\gb=\ga-\gg+\gd$. Similarly the equations of $[L_{\ga\gb}N_{\gg\gd}]$ can be
written
\begin{align*}
 & x_0-\ge^\ga x_1+\ge^{\ga+\gd}x_2-\ge^\gg x_3=f_{\ga+4,\ga+\gd,\gg+4}=0, \\
 & x_0^2+\ge^{\ga+\gd}x_0x_2+\ge^{2\ga+2\gd}x_2^2+\ge^{2\ga}x_1^2+\ge^{\ga+\gg}x_1x_3
+\ge^{2\gg}x_3^2=0.
\end{align*}
Thus the intersection of $X$ with $f_{\ga+4,\ga+\gd,\gg+4}=0$ is
\[ L_{\ga\gb}+N_{\gg\gd}+[L_{\ga\gb}N_{\gg\gd}] \]
where $\gb=\gg-\ga-\gd$.

We shall need some further intersections. If we write
\[ e'_\pm=x_0x_3\pm x_1x_2, \quad e''_\pm=x_0x_2\pm x_1x_3 \]
then the intersections of $X$ with the corresponding quadrics are as follows:
\begin{align*}
e'_-=0: & \quad L_{11}+L_{33}+L_{55}+L_{77}+M_{11}+M_{33}+M_{55}+M_{77}, \\
e'_+=0: & \quad L_{15}+L_{37}+L_{51}+L_{73}+M_{15}+M_{37}+M_{51}+M_{73}, \\
e''_-=0: & \quad L_{17}+L_{35}+L_{53}+L_{71}+N_{11}+N_{33}+N_{55}+N_{77}, \\
e''_+=0: & \quad L_{13}+L_{31}+L_{57}+L_{75}+N_{15}+N_{37}+N_{51}+N_{73}.
\end{align*}
Again, if we write
\begin{align*}
 & h'_{\ga\gb}=x_0^2-\ge^{2\ga}x_1^2-\ge^{2\gb}x_2^2+\ge^{2\ga+2\gb}x_3^2, \\
 & h''_{\ga\gb}=x_0^2-\ge^{2\ga}x_1^2-\ge^{2\gb}x_3^2+\ge^{2\ga+2\gb}x_2^2,
\end{align*}
which only depend on $\ga,\gb$ mod 4, then the intersection of $X$ with
$h'_{\ga\gb}=0$ is
\[ L_{\ga\ga}+L_{\ga,\ga+4}+L_{\ga+4,\ga}+L_{\ga+4,\ga+4}+
M_{\gb\gb}+M_{\gb,\gb+4}+M_{\gb+4,\gb}+M_{\gb+4,\gb+4} \]
and the intersection of $X$ with $h''_{\ga\gb}=0$ is
\[ L_{\ga,-\ga}+L_{\ga,4-\ga}+L_{\ga+4,-\ga}+L_{\ga+4,4-\ga}+
N_{\gb\gb}+N_{\gb,\gb+4}+N_{\gb+4,\gb}+N_{\gb+4,\gb+4}. \]
Moreover, on $X$ we have
\beq{E1} h'_{13}h'_{31}=-h'_{11}h'_{33}=2e'_+e'_-, \quad
h''_{13}h''_{31}=-h''_{11}h''_{33}=2e''_+e''_-. \end{equation}

There are a number of sets of 16 mutually skew curves of
genus 0 on $X$, of which a typical one is
\begin{align*}
 & M_{51},M_{33},M_{15},M_{77},\;[L_{33}M_{11}], [L_{33}M_{55}], [L_{15}M_{37}],
[L_{15}M_{73}], \\
 & N_{11},N_{37},N_{55},N_{73},\;[L_{57}N_{15}], [L_{57}N_{51}], [L_{71}N_{33}],
[L_{71}N_{77}].
\end{align*}
The map $X\rightarrow K$ which we shall exhibit identifies these 16 curves
with the 16 disjoint lines on $K=\Kum(A)$ that correspond to the points of
order 2 on $A$.

Now write
\[ D'=2L_{15}+M_{37}+M_{73}+[L_{31}N_{37}]+[L_{31}N_{73}] \]
and consider those principal divisors on $X$ of degree 16 which contain $D'$. Four
examples of them are given in the following table, which also names the associated
functions of $x_0,\ldots,x_3$ which give rise to them. The table is followed by
rational expressions for these functions; they can of course also be written as
polynomials, but the resulting formulae are unhelpful.
\begin{align*}
 & F_1: \quad D'+2L_{73}+M_{15}+M_{51}+[L_{57}N_{15}]+[L_{57}N_{51}], \\
 & F_2: \quad D'+2L_{55}+M_{33}+M_{77}+[L_{71}N_{33}]+[L_{71}N_{77}], \\
 & F_3: \quad D'+2L_{31}+N_{37}+N_{73}+[L_{15}M_{37}]+[L_{15}M_{73}], \\
 & F_4: \quad D'+2L_{17}+N_{11}+N_{55}+[L_{33}M_{11}]+[L_{33}M_{55}].
\end{align*}
Here we can take
\begin{align*}
 & F_1=\frac{f_{727}f_{763}f_{125}f_{161}e'_+(x_0-\ge x_1)(x_0-\ge^7x_1)}
{e''_+(x_2-\ge x_3)(x_2-\ge^7x_3)}, \\
 & F_2=\frac{f_{727}f_{763}f_{327}f_{363}h'_{13}(x_2-\ge^5x_3)}
{h''_{33}(x_2-\ge x_3)}, \\
 & F_3=f_{727}f_{763}f_{534}f_{570}, \\
 & F_4=\frac{f_{727}f_{763}f_{754}f_{710}e''_-h'_{13}(x_0-\ge x_1)(x_0-\ge^7x_1)}
{e'_-h''_{33}(x_2-\ge x_3)(x_2-\ge^3x_3)}.
\end{align*}
The divisors residual to $D'$ in the table of divisors associated with the $F_i$
have self-intersection 0 and are linearly equivalent, so they lie in a pencil,
which we denote by $\sP'$. Hence the restrictions of any three of
the $F_i/F_3$ to $X$ are linearly dependent. To find their linear dependence
relations, it is enough to consider for example their restrictions to $M_{13}$, and in
this way we find that on $X$
\[ F_3=\ge^2F_2-\ge(1+\ge^2)F_1, \quad F_4=-\ge(1+\ge^2)F_2+\ge^2F_1. \]

Let us also write
\[ D''=2L_{33}+M_{11}+M_{55}+[L_{75}N_{37}]+[L_{75}N_{73}] \]
and consider those principal divisors on $X$ of degree 16 which contain $D''$.
The table which follows, and the associated formulae, correspond to the ones
given above for the $F_i$.
\begin{align*}
 & G_1: \quad D''+2L_{11}+M_{33}+M_{77}+[L_{57}N_{15}]+[L_{57}N_{51}], \\
 & G_2: \quad D''+2L_{37}+M_{15}+M_{51}+[L_{71}N_{33}]+[L_{71}N_{77}], \\
 & G_3: \quad D''+2L_{75}+N_{37}+N_{73}+[L_{33}M_{11}]+[L_{33}M_{55}], \\
 & G_4: \quad D''+2L_{53}+N_{11}+N_{55}+[L_{15}M_{37}]+[L_{15}M_{73}].
\end{align*}
Here we can take
\begin{align*}
 & G_1=\frac{f_{367}f_{323}f_{125}f_{161}e'_-(x_0-\ge x_1)(x_0-\ge^3x_1)}
{e''_+(x_2-\ge^5x_3)(x_2-\ge^7x_3)}, \\
 & G_2=\frac{f_{367}f_{323}f_{327}f_{363}h'_{31}(x_0-\ge^3x_1)}
{h''_{33}(x_0-\ge^7x_1)}, \\
 & G_3=f_{367}f_{323}f_{754}f_{710}, \\
 & G_4=\frac{f_{367}f_{323}f_{534}f_{570}e''_-h'_{31}(x_2-\ge x_3)(x_2-\ge^3x_3)}
{e'_+h''_{33}(x_0-\ge x_1)(x_0-\ge^7x_1)}.
\end{align*}
The divisors residual to $D''$ in the table of divisors associated with the $G_i$
have self-intersection 0 and are linearly equivalent, so they lie in a pencil,
which we denote by $\sP''$. Hence the restrictions of any three of
the $G_i/G_3$ to $X$ are again linearly dependent. To find their linear dependence
relations, it is enough to consider for example their restrictions to $M_{13}$, and in
this way we find that on $X$
\[ G_3=\ge(1+\ge^2)G_1-\ge^2G_2, \quad G_4=-\ge(1+\ge^2)G_2+\ge^2G_1. \]

The restrictions to $X$ of $F_1G_2/F_2G_1$ and $F_3G_3/F_4G_4$ both have
divisors divisible by 2, so up to multiplication by a constant they are squares in $k(X)$.
Making use of (\ref{E1}) we find that
\[ \frac{F_1G_2}{F_2G_1}=2\left(\frac{e'_+}{h'_{13}}\right)^2,
\quad \frac{F_3G_3}{F_4G_4}=2\left(\frac{e''_+}{h''_{11}}\right)^2. \]
Thus if we write
\beq{E2} w_1=\frac{\ge}{1+\ge^2}\cdot\frac{F_2}{F_1}, \quad
w_2=\frac{\ge}{1+\ge^2}\cdot\frac{G_2}{G_1} \end{equation}
then we have
\begin{align*}
F_3/F_1=\ge(1+\ge^2)(w_1-1), & \quad F_4/F_1=-2\ge^2(w_1-\half), \\
G_3/G_1=-\ge(1+\ge^2)(w_2-1), & \quad G_4/G_1=-2\ge^2(w_2-\half).
\end{align*}
In particular we have a rational map $X\rightarrow K$ given by the equations (\ref{E2})
for $w_1,w_2$ together with
\beq{E3} z=\ge^3(1+\ge^2)\frac{e'_+}{h'_{13}}, \quad
          y=2\ge^2\frac{e''_+} {h''_{11 }}.
\end{equation}
But the curves $w_1$=constant and $w_2$=constant on $X$ are elements of $\sP',
\sP''$ respectively, so that their intersection has degree 4. In other words, if
$k$ contains $\ge$ then $[k(X):k(w_1,w_2)]=4$; and from this it follows that the
map $X\to K$ is actually birational. Since both $X$ and $K$ are minimal
models in their birational equivalence class, any birational map $X\to K$
is a biregular isomorphism. QED

\bigskip

\noindent Ecole Polytechnique F\'ed\'erale de Lausanne,
EPFL-SFB-IMB-CSAG, Station 8, CH-1015 Lausanne, Switzerland
\smallskip

\noindent Evis.Ieronymou@epfl.ch

\bigskip

\noindent Department of Mathematics, South Kensington Campus,
Imperial College London, SW7 2BZ England, U.K.

\smallskip

\noindent Institute for the Information Transmission Problems,
Russian Academy of Sciences, 19 Bolshoi Karetnyi,
Moscow, 127994 Russia
\medskip

\noindent a.skorobogatov@imperial.ac.uk

\bigskip

\noindent Department of Mathematics, Pennsylvania State University,
University Park, PA 16802, USA
\smallskip

\noindent Institute for Mathematical Problems in Biology, Russian
Academy of Sciences, Pushchino, Moscow Region, Russia

\medskip

\noindent zarhin@math.psu.edu

\end{document}